\RequirePackage{fix-cm}
\documentclass[smallextended]{svjour3}       
\smartqed  
\usepackage{graphicx}
\usepackage{amsmath}
\usepackage{amssymb} 

\begin{document}

\title{The Collatz Network
}
\subtitle{}


\author{Tobías Canavesi         
}


\institute{Physics Institute of La Plata \at
              Diagonal 113 63 y 64, CC67, 1900 La Plata, Argentina \\
              Tel.: +5491132531117\\
              ORCID: 0000-0003-2285-3103\\
              \email{tobiascanavesi@gmail.com}         
}

\date{Received: date / Accepted: date}

\maketitle

\begin{abstract}
Considering all possible paths that a natural number can take following the rules of the algorithm proposed in the Collatz conjecture we construct a graph that can be interpreted as an infinite network that contemplates all possible paths within the conjecture.  This allows us to understand why the minimal element of the Collatz orbit $x$ is equal to 1. Subsequently we define the extended Collatz conjecture equal to $o x+1$ when $x$ is odd and $x/2$ when $x$ is even with $o$ an odd number greater than $3$ and we show that there are infinite orbits in the extended Collatz conjecture that diverges. Finally, we find interesting theorems relating the Fibonacci sequence to prime numbers.
\end{abstract}

\section{Introduction}
\label{sec:introduction}
The Collatz problem, also known as the 3x+1 problem, is a very interesting an easy algorithm, and is one of those problems whose statement is very simple to understand.  Starting from a positive integer $x$, we begin to iterate in the following way, if it is even we divide it by two $x/2$ and if it is odd we multiply it by three and add one $3x+1$. The $3x + 1$ conjecture holds that, starting from any positive integer x, the repeated iteration of this function finally produces the value 1. We remit the reader to \cite{Lagarias} for an exhaustive discussion of this conjecture.\\
In section \ref{sec:paths} we study the conjecture from a different point of view to the one that normally is approached from number theory since we do it through a graph whose connections or edges are the possible paths that a natural number can take once the algorithm proposed in the conjecture is executed. \\
In section \ref{sec:Fibonacci} we consider all the possible paths that a natural number can take within the conjecture differentiating between even paths (resulting from dividing the number by two) and odd paths (resulting from multiplying the number by three and adding 1) and we show that in both cases the increase of even and odd paths results in the Fibonacci sequence. Finally we show the quotient between odd and even paths is equal to the inverse of the golden ratio.\\
In section \ref{sec:collatz} we show how the relationship found in section \ref{sec:Fibonacci} can be a great help in understanding the convergence observed in the Collatz conjecture. Finally we show that if we extend the conjecture by modifying it in such a way that when the number obtained is odd instead of multiplying by 3 we multiply by any odd number greater than 3, the conjecture has infinitely many divergent orbits.\\
In the last section \ref{sec:primes} we demonstrate interesting theorems relating prime numbers with the fibonacci sequence. 

\section{Paths in Collatz conjecture}
\label{sec:paths}
\begin{figure}[h]
	\centering
	\includegraphics[width=0.5\textwidth]{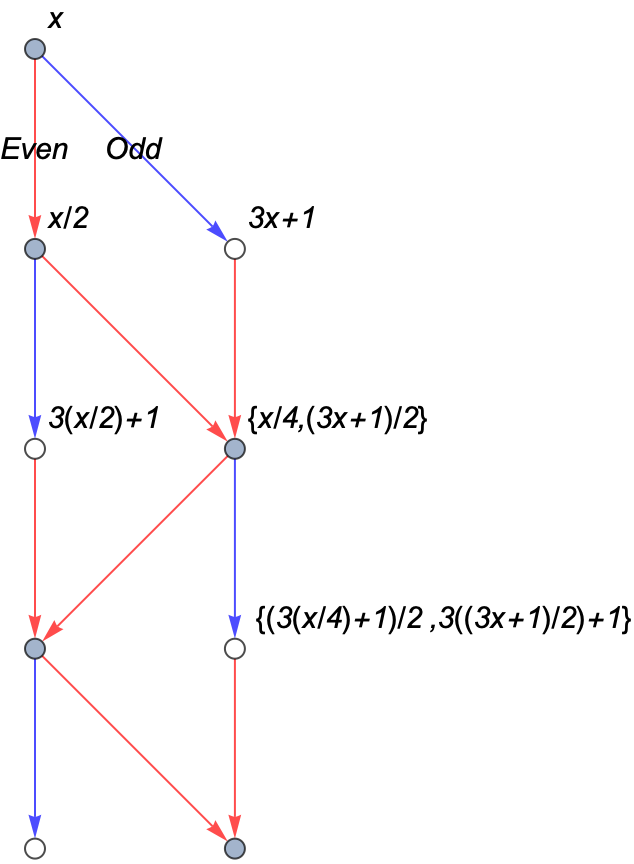}
	\caption{\label{fig:collatzgraph} This graph is a representation of the Collatz conjecture. Starting with $x\in \mathbf{N}$ number, the edges in blue represent a odd event and the red ones and even event. Each vertices of the graph represent one of the two possible operations in Collatz conjecture. The white vertices indicate that the result of that operation is a even number.}
\end{figure}

Starting from any natural number if it is even we divide it by two and if it is odd we multiply it by three and add one. The conjecture claims that starting from any natural number this algorithm always converges to 1. We can represent these mathematical operations at each node of a graph Fig.\hspace{0.1cm}\ref{fig:collatzgraph}, in the white nodes the arithmetic operation corresponds to take a number and multiply it by $3$ and adding $1$, and on the gray ones correspond to divide the number by $2$, the red edges of the graph occur when the result obtained at a node is even and the blue ones for the odd case. . The mathematically correct way of stating the Collatz conjecture is as follows \cite{Tao} \\

Define the Collatz map $\mathrm{Col} : \mathbb{N}+1 \to \mathbb{N}+1$ on the positive integers $\mathbb{N}+1 = \{1,2,3,\dots\}$ by setting $\mathrm{Col}(x)$ equal to $3x+1$ when $x$ is odd and $x/2$ when $x$ is even, and let $\mathrm{Col}_{\min}(x) := \inf_{x \in \mathbb{N}} \mathrm{Col}^n(x)$ denote the minimal element of the Collatz orbit $x, \mathrm{Col}(x), \mathrm{Col}^2(x), \dots$.\\
	The \textbf{Collatz conjecture} asserts that $\mathrm{Col}_{\min}(x)=1$ for all $x \in \mathbb{N}+1$. 
	
We can define different paths depending on the result obtained at each node, if the result is even we connect the following node with a red edge and if it is odd with a blue edge.

\begin{figure}[ht]
	\centering
	\includegraphics[width=0.2\textwidth]{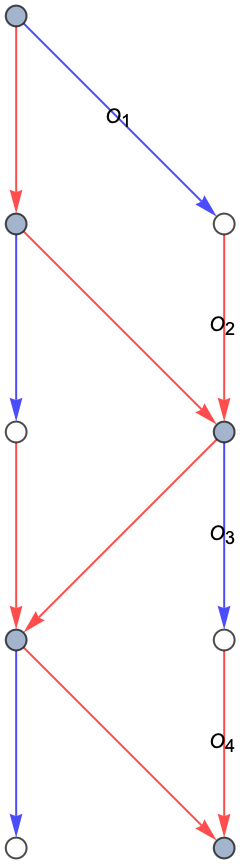}
	\caption{\label{fig:graphO} A uniform alternate color path starting from a odd number.}
\end{figure}

 We defined an \textit{uniform alternate path}, as a path in the graph whose edges are cyclic in the sense of a repetition of series of red and blue colors. We called recurrence of an alternate color path to the unit of recurrence, for example in the case of an RBRBRB... path the recurrence is RB. A example of this recurrence is in Fig.\hspace{0.1cm}\ref{fig:graphO}. If we consider an uniform alternate path with recurrence $RB$ we observe for this path that if we have $j$ red edges the maximum of blue edges is $j-1$ this is because always after a red edge comes a blue one and this path has the highest ratio between blue and red edges.
 
We defined an \textit{uniform path}, as a path in the graph whose edges are all from the same color.
There are two possible paths of this type the first is the path that comes from the powers of two that type of path that has all red edges is clearly ends at number $1$ as the last step in the iteration. The other possible uniform path corresponds to the case where all edges are blue, but it is easy to see that this is not possible since always after a blue edge comes a red one.

\begin{figure}[ht]
	\centering
	\includegraphics[width=0.2\textwidth]{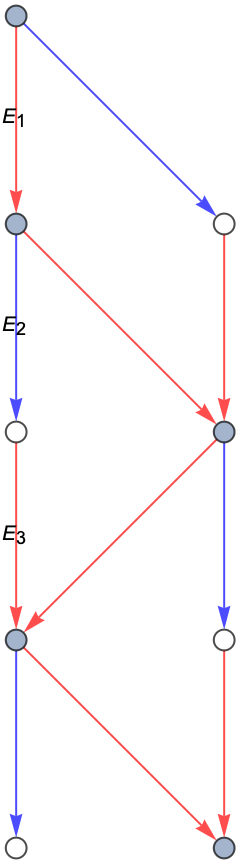}
	\caption{\label{fig:graphE} A uniform alternate color path starting from a even number.}
\end{figure}
Let focus our attention on the alternating path with $RB$ recurrence since this path as we said before is extreme in the sense that the ratio between red and blue edges is the highest possible. Starting with an arbitrary odd number $x\in\mathbf{N}$ we can make a general representation in terms of the result obtained at each node as in Fig.\hspace{0.1cm}\ref{fig:graphO}. We use the letter $P$ to keep track of how many times we apply the algorithm
\begin{equation}
x \Rightarrow v_1=P(3x+1),
\label{eq:1}
\end{equation}
the result is a even number, so be have to divide by two, so following the path.
\begin{equation}
P(3x+1)\Rightarrow v_2=P \frac{P(3x+1)}{2},
\label{eq:2}
\end{equation}
continuing with the path,
\begin{equation}
P\frac{P(3x+1)}{2} \Rightarrow v_3=P\left(3 P\frac{P(3x+1)}{2}  +1 \right).
\label{eq:3}
\end{equation}
we can characterize this for a vertex $j$ and rewrite these expressions in a general form:
\begin{equation}
\frac{P}{2}+\frac{3P^2}{2^2}+\frac{3^2P^3}{2^3}+\frac{3^3 P^4}{2^4}+...+\frac{3^{j-1}P^j}{2^{j}}(3 x +1).
\label{eq:general1}
\end{equation}
Where $j$, give the number of times that the red edge was applied, putting $P=1$.
\begin{equation}
\frac{1}{2}+\frac{3}{2^2}+\frac{3^2}{2^3}+\frac{3^3}{2^4}+...+\frac{3^{j-1}}{2^{j}}(3 x+1).
\label{eq:general2}
\end{equation}
So if we call $S$ to the sum of Eq.\ref{eq:general2},
\begin{equation}
S=\sum_{i=1}^{j-1} \frac{3^{i-1}}{2^{i}}+\frac{3^{j-1}}{2^{j}}(3x+1).
\label{eq:general3}
\end{equation}
Finally we find a general expression for this type of path. We get the same result starting from an even number, changing $3x+1\iff x$ and using the path in Fig.\hspace{0.1cm}\ref{fig:graphE}. So we find a general expression to represent this specific type of path which is going to be extremely important in Sec.\hspace{0.1cm}\ref{sec:collatz}.\\

\section{Fibonacci sequence in Collatz conjecture}\label{sec:Fibonacci}
\begin{figure}[h]
	\centering
	\includegraphics[width=1\textwidth]{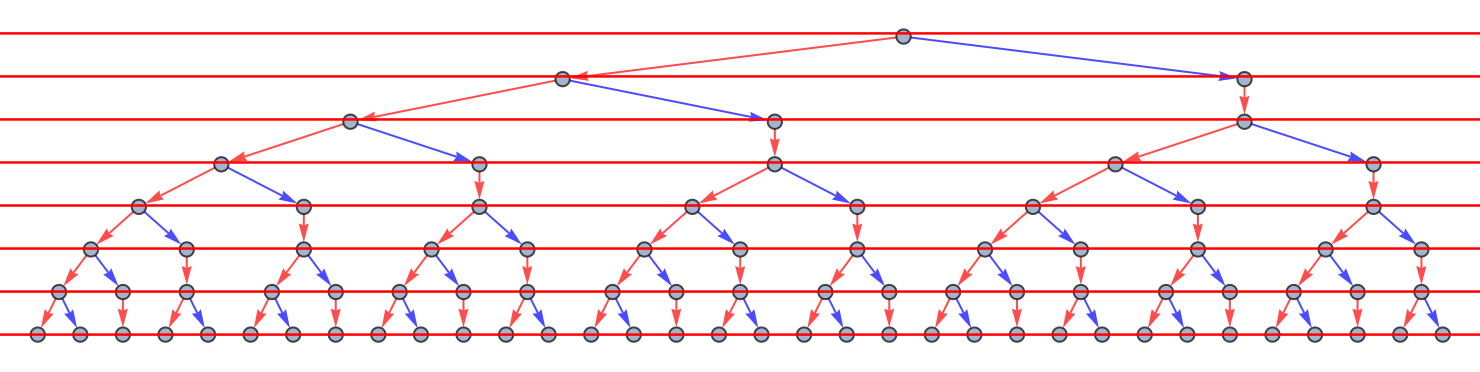}
	\caption{\label{fig:fibonacci} In this graph we can observe all the possible branches of evolution(iteration) of the conjecture of Collatz for 7 iterations. The horizontal red lines indicate the evolution from zero or no iteration to iteration 7.}
\end{figure}

Fibonacci numbers are a fascinating sequence of integers numbers discovered by the mathematician Leonardo Fibonacci, related to the shapes of flower petals, tree branches, rabbit birth rates, and many other natural phenomena.  
The Fibonacci sequence is one of the oldest known recursive sequence, which is a sequence where each successive term can only be found through performing operations on previous terms.
We can find this famous sequence also in the Collatz conjecture. For this we will interpret the conjecture as a graph which contains all the possible paths that can occur within it. The graph is constructed by considering all possible paths at each node, this means that from a red edge can continue a red or blue edge and from a blue edge can continue only a red one. So starting with a natural number $x \in \mathbb{N}+1$ we can study the evolution of the red and blue edges considering all possible paths Fig.\hspace{0.1cm}\ref{fig:fibonacci}, we found the evolution for the edges in Tab.\ref{tab:table1}.\\

\begin{table}[h!]
	\begin{center}
		\caption{Table indicating the number of red and blue edges per iteration. Where Fn indicates the }
		\label{tab:table1}
		\begin{tabular}{c|c|c|c}
			\hline
			\textbf{Nº of iterations} & \textbf{Nº of red edges} & \textbf{Nº of blue edges} &\textbf{Quotient B/R}\\
			\hline
			0 & 0 & 0&-\\
			\hline
			1 & 1 & 1& 1\\
			\hline
			2 & 2 & 1& $1/2$ \\
			\hline
			3 & 3 & 2& $2/3$ \\
			\hline
			4 & 5 & 3& $3/5$ \\
			\hline
			5 & 8 & 5& $5/8$ \\
			\hline
			6 & 13 & 8& $8/13$ \\
			\hline
			7 & 21 & 13& $13/21$ \\
			... & ... & ...& ...\\
			$\infty$ & ... & ...& $1/\varphi$\\
			\hline
		\end{tabular}
	\end{center}
\end{table}
We find therefore the Fibonacci sequence in the growth of red and the blue edges, but out of phase in one iteration. If we define $R(n)$ as the number of red edges and $B(n)$ the number of blue edges in the $n$ iteration. We find that 
\begin{eqnarray}
\label{eq:goldenratio}
\lim_{n \to +\infty} \frac{B(n)}{R(n)}=\frac{1}{2} \left(\sqrt{5}-1\right)=\frac{1}{2} \left(1+\sqrt{5}\right)-1
\end{eqnarray}
Which is nothing but the golden ratio minus one. This is indicating that we can find a linear relationship in the case of ${n \to +\infty}$ with slope $m=\left(\frac{1}{2} \left(1+\sqrt{5}\right)-1\right)$. We can do some simple algebraic computations in Eq.\hspace{0.1cm}\ref{eq:goldenratio} for this limit case
\begin{eqnarray}
\label{eq:goldenratiolineal}
B(n)=\left(\frac{1}{2} \left(1+\sqrt{5}\right)-1\right)R(n)=\frac{R(n)}{\varphi}=\\
R(n)(-1+\varphi)=R(n)\varphi -R(n)
\end{eqnarray}
so
\begin{eqnarray}
\label{eq:goldenratiolineal2}
B(n)+R(n)=R(n)\varphi 
\end{eqnarray}
So we find the Fibonacci sequence in this interpretation of the conjecture.\\
It is important to understand that taking the limit $n\rightarrow +\infty$ means that in this case the result will be an infinite graph that includes all possible paths within the conjecture.

\section{Collatz conjecture}\label{sec:collatz}
We can write a general expression to take into account the relation between the grow of blue and red edges. Considering that the blue edges represent the powers of 3 and the red ones those of 2. 
\begin{equation}\label{eq:coc}
\frac{3^{F_{n}}}{2^{F_{n}+F_{n-1}}}=\frac{3^{F_{n}}}{2^{F_{n+1}}}=a_{n}.
\end{equation}
\begin{figure}[h]
	\centering
	\includegraphics[width=0.8\textwidth]{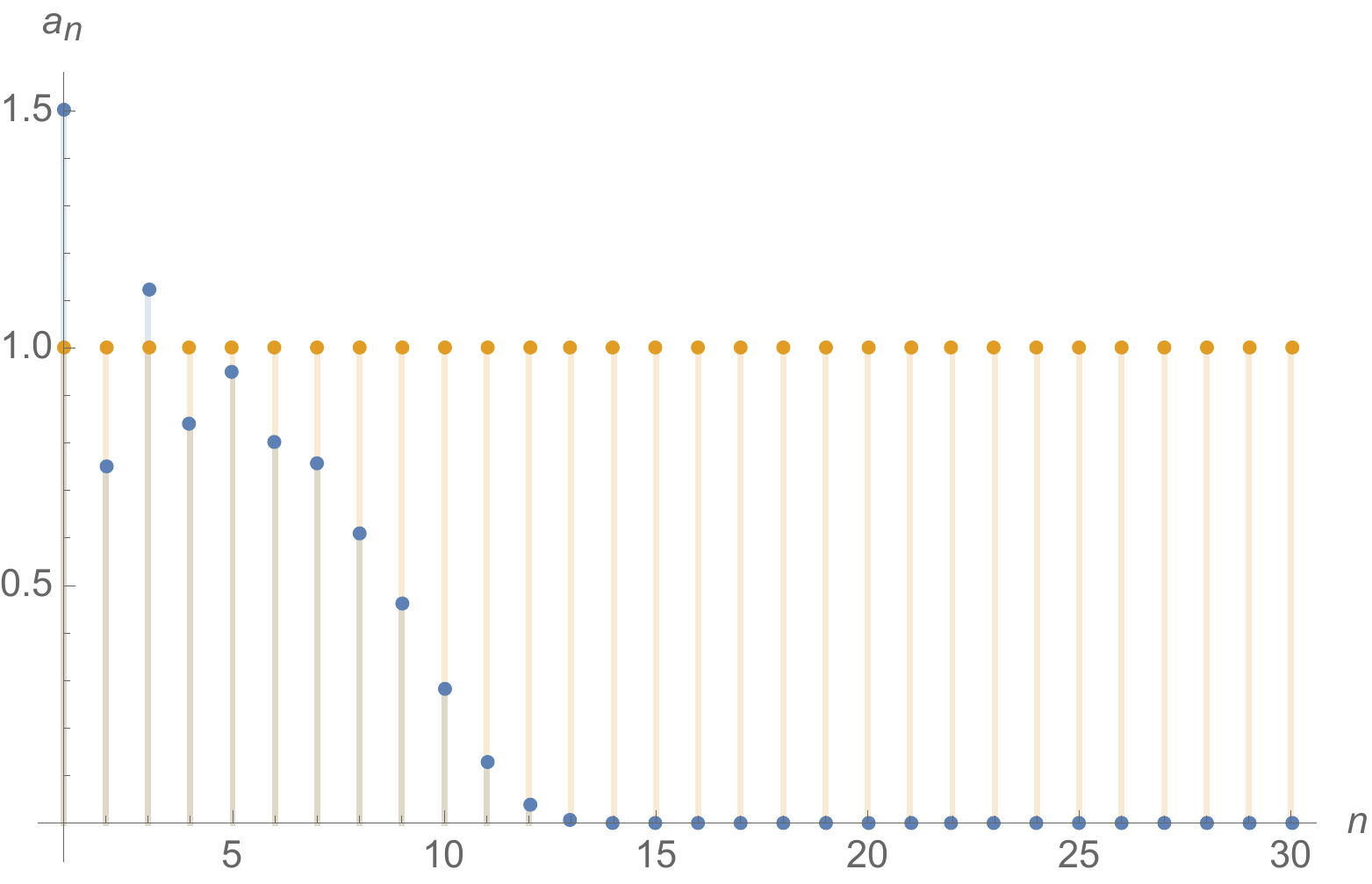}
	\caption{\label{fig:succ} A discrete plot of $a_n$ (blue points) and in the yellow points represent the line $y=1$. }
\end{figure}
We can see how the succession evolves for different $n$ in Fig.\hspace{0.1cm}\ref{fig:succ}.
\begin{lemma}
The succession $\frac{3^{F_{n}}}{2^{F_{n+1}}}=a_{n}$ converge to $0$.
\end{lemma}
\textrm{Proof: For this let prove that $a_{n}<b_{n}=1/n$.}
\begin{eqnarray}
\frac{3^{F_{n}}}{2^{F_{n+1}}}\stackrel{?}{<}\frac{1}{n}\\
F_{n}\log(3)-F_{n+1}\log(2)\stackrel{?}{<}-\log(n)\\
F_{n}\log(3)+\log(n)\stackrel{?}{<}F_{n+1}\log(2)\\
\frac{F_{n}}{F_{n+1}}\frac{\log(3)}{\log(2)}+\frac{\log(n)}{F_{n+1}\log(2)}\stackrel{?}{<}1
\end{eqnarray}
\textrm{taking the limit $n\rightarrow \infty$ then}
\begin{eqnarray}
\frac{1}{\varphi}\frac{\log(3)}{\log(2)}+0\stackrel{?}{<}1\\
\frac{1}{\varphi}\stackrel{?}{<}\frac{\log(2)}{\log(3)}\\
\frac{1}{\varphi}<\frac{\log(2)}{\log(3)}
\end{eqnarray}
\textrm{So the inequality is fulfilled and therefore the lemma is proven because is well know that $\lim_{n\rightarrow \infty} b_{n}=\lim_{n\rightarrow \infty}\frac{1}{n}=0$. $\qed$}\\

\begin{figure}[h]
	\centering
	\includegraphics[width=0.8\textwidth]{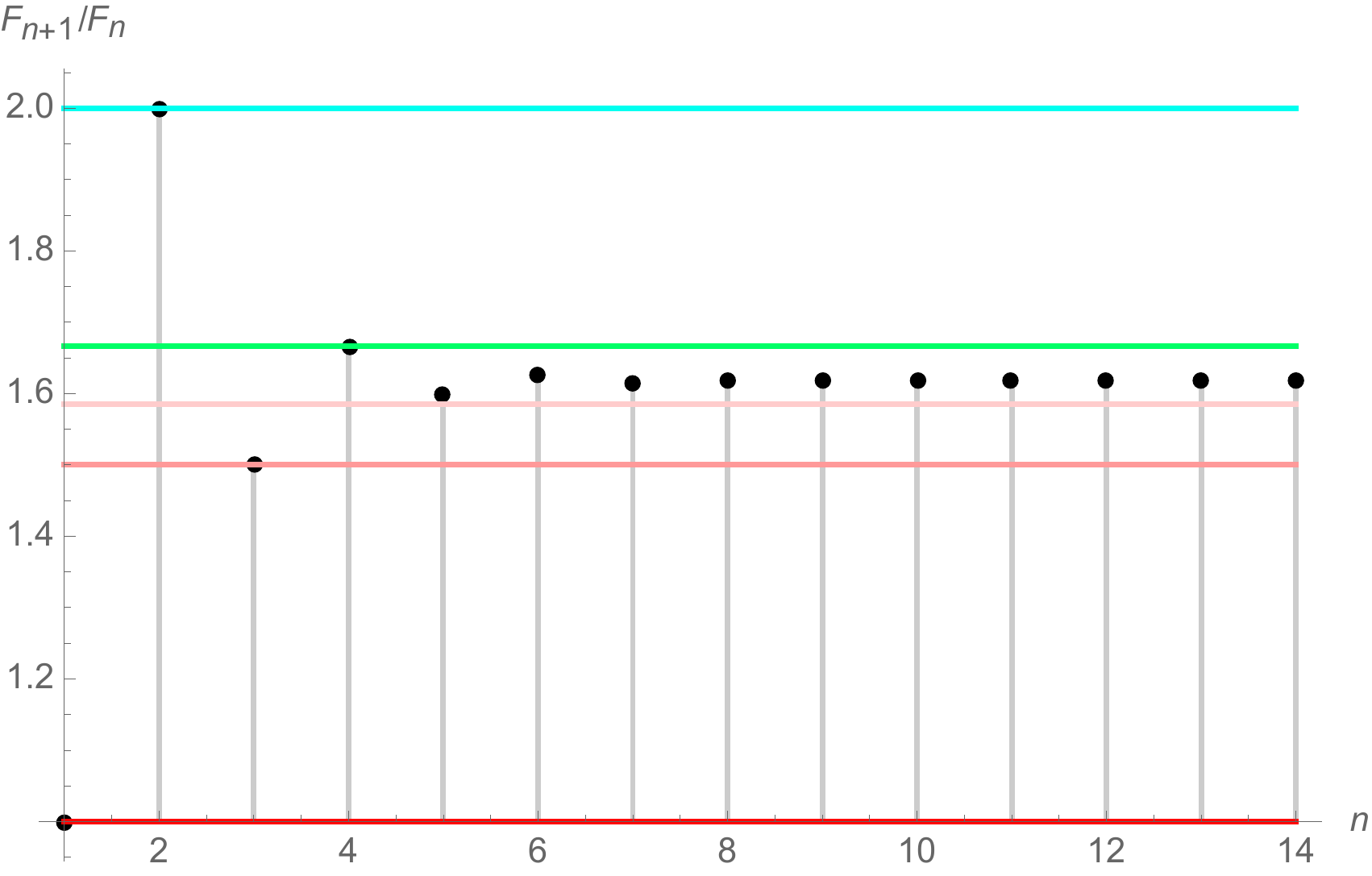}
	\caption{\label{fig:rb} A discrete plot of $\frac{F_{n+1}}{F_{n}}$ (black points) and and different horizontal lines given by the equations, $y=\frac{F_{1+1}}{F_{1}}$(red) $y=\frac{F_{2+1}}{F_{2}}$(lightblue), $y=\frac{F_{3+1}}{F_{3}}$(lightred), $y=\frac{F_{4+1}}{F_{4}}$(lightgreen) and $y=\frac{\log(3)}{\log(2)}$(pink). }
\end{figure}

\begin{theorem}
$\forall x \in \mathbb{N}+1$ all the possible paths in Collatz conjecture converge to 1 or $\mathrm{Col}_{\min}(x)=1$ for all $x \in \mathbb{N}+1$.
\end{theorem}
\textrm{Proof:} As we explained in the previous section all the paths in the Collatz conjecture are contained within the graph Fig.\hspace{0.1cm}\ref{fig:fibonacci} although this also contains paths that include real numbers and that are not part of the paths in the conjecture this does not generate any problem as we will show later.  Considering all possible paths starting with some $x \in \mathbb{N}+1$ that follows the arithmetic operations consider in the Collatz conjecture. In Sec. \ref{sec:Fibonacci} we find a relation for the grow of blue and red edges for the case $n\rightarrow +\infty$ Eq.\hspace{0.1cm}\ref{eq:goldenratiolineal}
	\begin{equation}
	B(n)=\frac{R(n)}{\varphi},
	\end{equation}	
	where $\varphi$ is the golden ratio. Let's start by analyzing how the ratio between blue and red edges evolves, 
\begin{eqnarray*}
	j=1\hspace{1cm}\frac{1}{\varphi_{1}}=\frac{F_{1}}{F_{2}}=1\\
	j=2\hspace{1cm}\frac{1}{\varphi_{2}}=\frac{F_{2}}{F_{3}}=\frac{1}{2}\\
	j=3\hspace{1cm}\frac{1}{\varphi_{3}}=\frac{F_{3}}{F_{4}}=\frac{2}{3}\\
	...\\
	\lim\limits_{n\rightarrow \infty}\frac{F_{j}}{F_{j+1}}=\frac{1}{\varphi}\\
\end{eqnarray*}
where $1/\varphi=\lim\limits_{n\rightarrow \infty} 1/\varphi_{n}=\lim\limits_{n\rightarrow \infty} \frac{F_{n}}{F_{n+1}}$. The minimum value that could take $1/\varphi_{n}$ is $1/\varphi_{2}$ Fig.\hspace{0.1cm}\ref{fig:rb}
\begin{eqnarray}\label{eq:ineq}
\varphi_n < \varphi_2 \\
2^{\varphi_n} < 2^{\varphi_2}\\
2^{\varphi_n} -3< 2^{\varphi_2}-3\\
\frac{1}{2^{\varphi_n} -3}>\frac{1}{2^{\varphi_2}-3}
\end{eqnarray}
then $\forall n \in \mathbb{N}+1$ then $\varphi_{2} > \varphi_{n}$ therefore the minimum for $\frac{1}{\varphi_n}$ is when $n=2$.

\begin{figure}[h]
	\centering
	\includegraphics[width=0.8\textwidth]{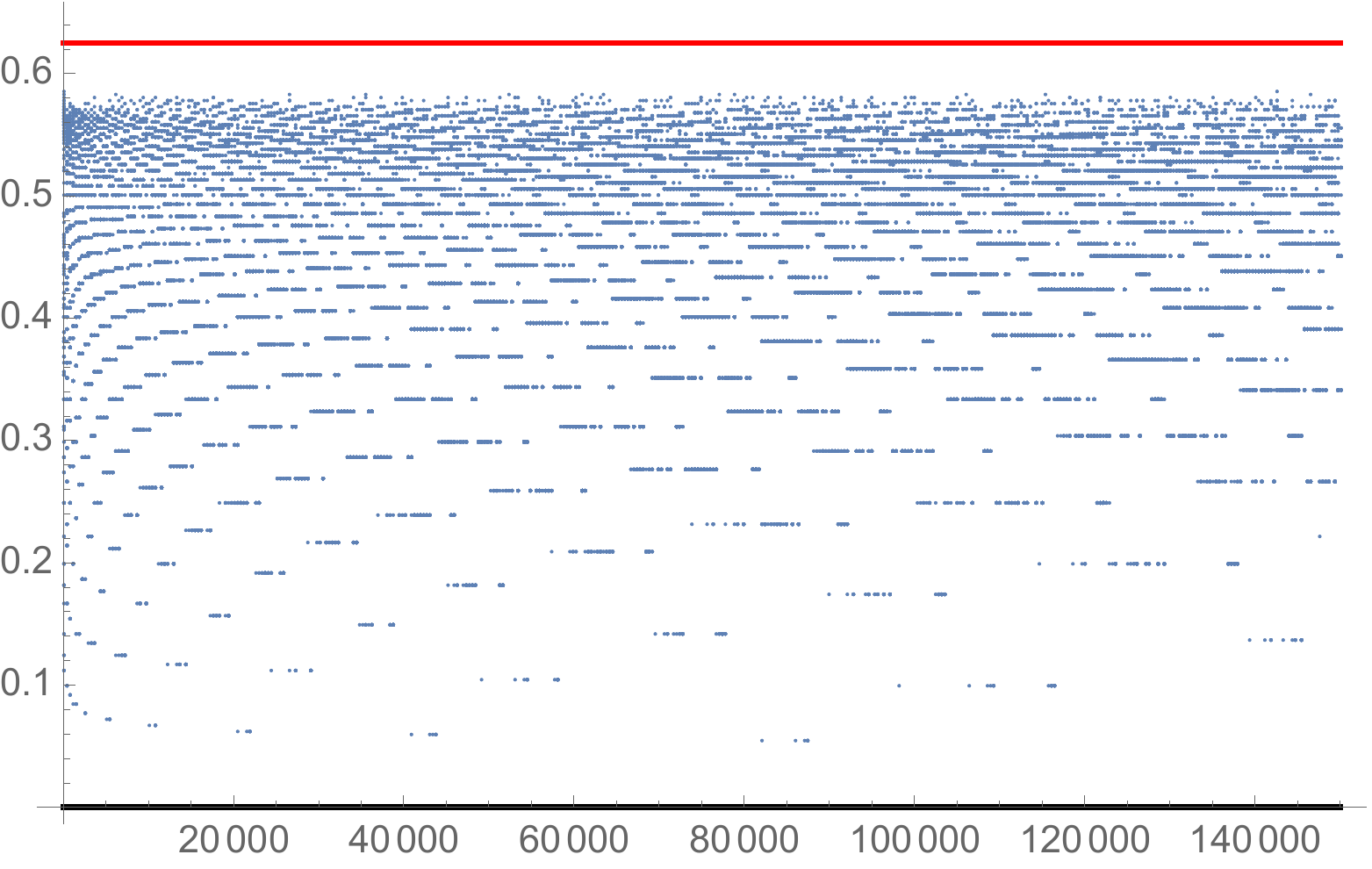}
	\caption{\label{fig:quotient} Ratio between the number of blue edges and the number of red edges ($y=B/R$) for the first $150000$ natural numbers in the Collatz conjecture. The red line corresponds to $y=\frac{5}{8}$ and the black line to $y=0$}
\end{figure}

The next step in the construction of our proof is to ask when a succession of the form $\frac{3^{n}}{2^{k n}}=\left(\frac{3}{2^{k}}\right)^{n}$ converges, this happen for $k> \frac{\log(3)}{\log(2)}$, therefore  $F_{n+1}/F_{n}>k$ except for the finite cases $n=1,3$. As any path in the conjecture can be written in a general form as  Eq.\hspace{0.1cm}\ref{eq:general3} we can perform this sum over all the possible paths for each $n$ and use the relation between red and blue edges given by 
 Eq.\hspace{0.1cm}\ref{eq:goldenratio}. Note that $j$ now sums over all paths in each iteration, not over a single particular path
\begin{equation}\label{eq:sumallpaths}
S=\sum_{i=1}^{j-1}\frac{3^{i-1}}{2^{i\varphi_n}}+\frac{3^{j-1}}{2^{j\varphi_n}}x=\frac{1}{3} 2^{-\varphi_n j} \left(3^j x+\frac{3\ 2^{j \varphi_n}-3^j 2^{\varphi_n }}{2^{\varphi_n }-3}\right),
\end{equation}
for the case starting with a odd we have to change $x\Longleftrightarrow 3x+1$. We can take this limit since there are infinite natural numbers from which I can begin to iterate on the conjecture so there are infinitely many possible paths. Therefore the limit $j\rightarrow \infty$ in Eq.\ref{eq:sumallpaths}
\begin{equation}
\label{eq:limitfibo4}
\lim\limits_{j\rightarrow \infty} S=\lim\limits_{j\rightarrow \infty}  \frac{1}{3} 2^{-\varphi_n j} \left(3^j x+\frac{3\ 2^{j \varphi_n
	}-3^j 2^{\varphi_n }}{2^{\varphi_n }-3}\right)=\frac{1}{2^{\varphi_n}-3}
\end{equation}
so that the infinite paths have the minimum in $n=2$ so $\frac{F_{2+1}}{F_{2}}=\varphi_2$ see Eq.\ref{eq:ineq}. Therefore changing $\varphi_n \Longleftrightarrow \varphi_{2}$ in Eq.\ref{eq:limitfibo4}
\begin{equation}
\label{eq:limitfibo5}
\frac{1}{2^{\varphi_2}-3}=\frac{1}{2^{2}-3}=1
\end{equation}
then for all the paths in the Collatz conjecture $\mathrm{Col}_{\min}(x)=1$ for all $x \in \mathbb{N}+1$. $\qed$\\
\begin{corollary}
There is a cut off for the quotient of blue and red edges given by $max\{\frac{B(n)}{R(n)}\}=\frac{5}{8}$ $\forall n$ except $n=1$ and $n=3$ \\
We can justify the exclusion of $n=1$ and $n=3$ in two ways, the first is by taking into account that they would give unbound orbits, but we know from \cite{Tao} that this is not possible. And the second is noticing that there are no paths in the Collatz conjecture that have such a relationship between blue and red edges.\\
	The maximum value for $\frac{F_{n}}{F_{n+1}}$ if we exclude $n=1$ and $n=3$ it is for the case $n=5$, $\frac{F_{5}}{F_{5+1}}=\frac{5}{8}$, so for the theorem above we know that this is a maximum cut off for the possible quotient between red and blue edges. In Fig.\hspace{0.1cm}\ref{fig:quotient} we can see the quotient between blue and red edges for the first $150000$ natural numbers in Collatz conjecture. $\qed$\\
\end{corollary}

Extending the Collatz conjecture to the form $ox+1$ where $o$ are all the odd numbers greater than 3, we can rewrite its definition as \\

Define the extended Collatz map $\mathrm{Col_x} : \mathbb{N}+1 \to \mathbb{N}+1$ on the positive integers $\mathbb{N}+1 = \{1,2,3,\dots\}$ by setting $\mathrm{Col_{x}}(x)$ equal to $o x+1$ when $x$ is odd and $x/2$ when $x$ is even, and let $\mathrm{Col_{x}}_{\min}(x) := \inf_{x \in \mathbb{N}} \mathrm{Col_x}^n(x)$ denote the minimal element of the Collatz orbit $x, \mathrm{Col_x}(x), \mathrm{Col_x}^2(x), \dots$.\\
\begin{theorem}
There are infinite orbits in the extended Collatz conjecture that diverges. 
\end{theorem}
\textrm{Proof:} We can do the same reasoning that the theorem above, but in this case $k\geq \frac{\log(o)}{\log(2)}$, if we called $O$ the set of all the odd numbers greater than 3 the $min\{O\}=5$ so we find an analog of Eq.\ref{eq:limitfibo4} but in this case we need $\varphi'>\frac{\log(5)}{\log(2)}$ for the series to converge to
	\begin{equation}
	\label{eq:limitfibo6}
    \frac{1}{2^{\varphi'}-5}
	\end{equation}
	but $ \frac{\log(5)}{\log(2)}>2=max\{\varphi_{n}=F_{n+1}/F_{n}\}$ so $\varphi'> F_{n+1}/F_{n}$ $\forall n\in \mathbb{N}+1$ so we notice that not all the paths converge, there are infinite paths that diverge. We can extend this for all $o>3$ $\qed$.

\section{Primes in Collatz conjecture}\label{sec:primes}
And interesting thing to study is the possible number of combinations of red and blue edges in each iteration
\begin{equation}
\binom{R(n)+B(n)}{B(n)}=\binom{F_{n+1}+F_{n}}{F_{n}}=\binom{F_{n+2}}{F_{n}}.
\end{equation}
This leads us to very relevant theorems:
\begin{theorem}\label{theorem1}
Let $A_{n+2}$ be the set of all the primes $\leq F_{n+2}$ (the equal holds when $F_{n+2}$ is prime) then $F_{n+2}!$ is divisible for all $p\in A_{n+2}$ for $n\geq 1$.
\end{theorem}
\textrm{Proof:} As $F_{n+2}!$ is simply $\prod_{q=2}^{F_{n+2}} q$  considering the fundamental theorem of arithmetic we know that there is a single prime factorization for each element of the product, so all the primes between $2 \leq p \leq F_{n+2} $ will be in $A_{n+2}$ and since if a prime $p$ is a divisor of $F_{n+2}!$ then $p$ must be a divisor of at least one number between $1,2,...F_{n+2}$(Euclid's Lemma) so $\forall p \in A_{n+2}$, $\max\{A_{n+2}\}\leq F_{n+2}$  we conclude that $\forall p \in A_{n+2}$ $p$ divide $F_{n+2}!$. $\qed$  \\
\begin{theorem}\label{theorem2}
 $F_{n+1}!$ has at least one divisor more than $F_{n}!$. Let $A_{n+1}$ be the set of all the primes $\leq F_{n+1}$ and $A_{n}$ be the set of all the primes $\leq F_{n}$ then $A_{n+1}>A_{n}$.
 \end{theorem}
\textrm{Proof:} From \cite{Boase} we know that if $n\geq 2$, then with the exception of $n = 5$ and $n = 11$, $F_{n+1}$ is divisible by some prime $p$ which does not divide any $F_{n}$ and $n<n+1$ this implies that $A_{n+1}$ has at least one element more than $A_{n}$. $\qed$
\begin{corollary}
We can build all the prime numbers through the prime factorization of the factorial of Fibonacci numbers.
\end{corollary}
\textrm{Proof:} From theorem \ref{theorem1} we know that $F_{n+2}!$ is divisible $\forall p \in A_{n+2}$ where $A_{n+2}$ contains all the prime numbers $\leq F_{n+2}$. For understand this let start with $n=1$ so $F_{3}=2$ and $A_{3}={2}$, for $n=2$ then $F_{4}=3$ and $A_{4}={3}$,for $n=3$ then $F_{5}=5$ and $A_{5}={5}$ we can continue with this process, so if we define $A:=A_{3}\cup A_{4}\cup ...\cup A_{n+2}$ this will be the set of all the primes $\leq F_{n+2}$ $\qed$.
\begin{theorem}\label{theoremfibo}
Let be $F_{n}$ the Fibonacci sequence and $\binom{F_{n+2}}{F_{n}}$ the possible number of combinations of red and blue edges in each iteration starting with $n=1$, the prime factorization of $\binom{F_{n+2}}{F_{n}}$ gives all the prime numbers ordered in an increasing way. Another way of stating this theorem would be, let $A'_{n+2}$ be the set of all the primes between $F_{n+1}$ and $F_{n+2}$ ( $\max  \{A'_{n+2} \}=F_{n+2}$ when $F_{n+2}$ is prime ) we can construct all the primes starting with $n=1$ and continuing $n=2,3,...,$ as a $A'_{1+2}\cup A'_{2+2}\cup...A'_{n+2}$. And the last and the more precise statement would be: Let $A'_{n+2}$ be a set of primes between the Fibonacci sequence $F_{n+1}$ and $F_{n+2}$ such that $A'_{n+2}:=\{p_{1},p_{2},...,p_{p}\}$ with $p_{p}=F_{n+2}$ if $F_{n+2}$ is prime and $p_{1}<p_{2}<...<p_{p}$, then $\binom{F_{n+2}}{F_{n}}$ is divisible $\forall p \in A'_{n+2}$ and $n\geq 1$.
\end{theorem}
\newpage
\textbf{Preliminary lemmas:}
\begin{lemma}
Let n be a positive integer. Then 
\begin{equation}
gcd \left(\{ \binom{n}{k} : 1 \leq k\leq n, gcd(k,n)=1 \} \right)=n.	
\end{equation}
\textrm{Proof:} Following \cite{Siao}, let $S_{n}$ be 
\begin{equation}
S_{n}=\{ \binom{n}{k} : 1 \leq k\leq n, gcd(k,n)=1 \}. 
\end{equation}
So $n\in S_{n}$ and so $gcd(S_{n})|n$. Let $k$ be an integer with $1\leq k\leq n$ and $gcd(k,n)=1$ and multiplying $k$ with the binomial $\binom{n}{k}$
\begin{equation}
k\binom{n}{k}=k\frac{n!}{k!(n-k)!}=\frac{n!}{(k-1)!(n-k)!}=n\frac{(n-1)!}{(k-1)!(n-k)!}=n\binom{n-1}{k-1}
\end{equation}
So $n$ divides $k\binom{n}{k}$. But $gcd(n,k)=1$ so n is coprime with k. So n divides $\binom{n}{k}$ then $n|gcd(S_{n})$ then
\begin{equation}
gcd \left( S_{n}\right)=n.  \qed
\end{equation}
\end{lemma}

\begin{lemma}\label{lemma3}
Let $F_{n}$ be the Fibonacci sequence,
\begin{eqnarray}
gcd \left(\{ \binom{F_{n+2}}{F_{n}} : 1 \leq F_{n}\leq F_{n+2},\hspace{0.1cm} gcd(F_{n+2},F_{n})=1 \} \right)=F_{n+2}.
\end{eqnarray}
\end{lemma}
\textrm{Proof:} Let start from showing that $gcd(F_{n+2},F_{n})=1$
\begin{eqnarray}
gcd(F_{n+2},F_{n})=gcd(F_{n+1}+F_{n},F_{n})=gcd(F_{n+1},F_{n})=\\
gcd(F_{n-1}+F_{n},F_{n})=gcd(F_{n-1},F_{n})
\end{eqnarray}
in each equality we used $gcd(i,i k+j)=gcd(i,j)$. This is true for all $n\geq 1$. Since $gcd(F_{1-1},F_{1})=gcd(0,1)=1$, this implies that $gcd(F_{n+2},F_{n})=1$ for all $n\geq 1$. Let $\mathbb{F}_{n}$ be
\begin{equation}
\mathbb{F}_{n}=\{ \binom{F_{n+2}}{F_{n}} : 1 \leq F_{n}\leq F_{n+2}, \hspace{0.1cm}gcd(F_{n+2},F_{n})=1 \} 
\end{equation}
So using an the same reasoning from lemma 1, $F_{n+2}\in \mathbb{F}_{n}$ and so $gcd(\mathbb{F}_{n})|F_{n+2}$. Let $F_{n}$ be an integer with $1\leq F_{n}\leq F_{n+2}$ and $gcd(F_{n+2},F_{n})=1$ and multiplying by $F_{n}$ the binomial $\binom{F_{n+2}}{F_{n}}$
\begin{eqnarray}
F_{n}\binom{F_{n+2}}{F_{n}}=F_{n}\frac{F_{n+2}!}{F_{n}!(F_{n+2}-F_{n})!}=\frac{F_{n+2}!}{(F_{n}-1)!(F_{n+2}-F_{n})!}=\\
F_{n+2}\frac{(F_{n+2}-1)!}{(F_{n}-1)!(F_{n+2}-F_{n})!}=F_{n+2}\binom{F_{n+2}-1}{F_{n}-1}
\end{eqnarray}
So $F_{n+2}$ divides $F_{n}\binom{F_{n+2}}{F_{n}}$. But $gcd(F_{n+2},F_{n})=1$ so $F_{n+2}$ is coprime with $F_{n}$. So $F_{n+2}$ divides $\binom{F_{n+2}}{F_{n}}$ then $F_{n+2}|gcd(\mathbb{F}_{n})$ then
\begin{equation}
gcd \left( \mathbb{F}_{n}\right)=F_{n+2}. \qed
\end{equation}
\begin{corollary}
Let $F_{n}$ be the Fibonacci sequence,
\begin{eqnarray}
gcd \left(\{ \binom{F_{n+1}}{F_{n}} : 1 \leq F_{n}\leq F_{n+1},\hspace{0.1cm} gcd(F_{n+1},F_{n})=1 \} \right)=F_{n+1}.
\end{eqnarray}
\end{corollary}
\textrm{Proof:} As the condition $gcd(F_{n+1},F_{n})=1$ is fulfilled the proof is analog to the lemma 2. \\
	
\textrm{Theorem \ref{theoremfibo} proof:}
Let's start by looking at some particular cases.\\
If $n=1$, $\binom{F_{3}}{F_{1}}$, $F_{3}=2$ and $F_{1}=1$ then $A_{3}=\{2\}$.\\
If $n=2$, $\binom{F_{4}}{F_{2}}$, $F_{4}=3$ and $F_{2}=1$ then $A_{4}=\{3\}$.\\
If $n=3$, $\binom{F_{5}}{F_{3}}$, $F_{5}=5$ and $F_{3}=2$ then $A_{5}=\{5\}$.\\
If $n=4$, $\binom{F_{6}}{F_{4}}$, $F_{6}=8$ and $F_{4}=3$ then $A_{6}=\{7\}$. \\
If $n=5$, $\binom{F_{7}}{F_{5}}$, $F_{7}=13$ and $F_{5}=5$ then $A_{6}=\{11,13\}$. \\
If $n=6$, $\binom{F_{8}}{F_{6}}$, $F_{8}=21$ and $F_{6}=8$ then $A_{7}=\{17,19\}$. \\
We can generalized this in the following form, we know from theorem \ref{theorem1}, that $F_{n+2}!$ is divided by all the primes $\leq F_{n+2}$ so  $F_{n+1}!$ is divided by all the primes $\leq F_{n+1}$ and  $F_{n}!$ is divided by all the primes $\leq F_{n}$ and we probe in lemma \ref{lemma3} that $gcd(F_{n+2},F_{n})=gcd(F_{n+2},F_{n+1})=gcd(F_{n+1},F_{n})=1$ if we look to the binomial coefficient 
\begin{eqnarray}
\binom{F_{n+2}}{F_{n}}=\frac{F_{n+2}!}{F_{n}!(F_{n+2}-F_{n})!}=\frac{F_{n+2}!}{F_{n}!F_{n+1}!}
\end{eqnarray}
let $A_{n}:= \{ p_{1},p_{2},...,p_{q} \}$ the prime divisors of $F_{n}!$, $A_{n+1}:= \{ p_{1},p_{2},...,p_{r} \}$ the prime divisors of $F_{n+1}!$ and $A_{n+2}:= \{ p_{1},p_{2},...,p_{s} \}$ the prime divisors of $F_{n+2}!$. We know from theorem \ref{theorem1} and theorem \ref{theorem2} that $A_{n} < A_{n+1} < A_{n+2} $ so this implies that there are primes between $F_{n+1}$ and $F_{n+2}$. We can write $F_{n+2}$, $F_{n+1}$ and $F_{n}$ in their respective prime factorization.
\begin{eqnarray*}
F_{n}!=p^{\alpha_{1}}_{1}p^{\alpha''_{2}}_{2}...p^{\alpha''_{q}}_{q} \hspace{0.3cm}\text{with}\hspace{0.3cm} p_{q} \leq F_{n} \\
F_{n+1}!=p^{\alpha'_{1}}_{1}p^{\alpha'_{2}}_{2}...p^{\alpha'_{q}}...p^{\alpha'_{r}}_{r} \hspace{0.3cm}\text{with}\hspace{0.3cm} p_{r} \leq F_{n+1} \\	
F_{n+2}!=p^{\alpha_{1}}_{1}p^{\alpha_{2}}_{2}...p^{\alpha_{q}}_{q}...p^{\alpha_{r}}_{r}...p^{\alpha_{s}}_{s} \hspace{0.3cm}\text{with}\hspace{0.3cm} p_{s} \leq F_{n+2} \\
\end{eqnarray*}
so the binomial take the form
\begin{eqnarray*}
\binom{F_{n+2}}{F_{n}}=\frac{F_{n+2}!}{F_{n}!F_{n+1}!}=\\
=\frac{p^{\alpha_{1}}_{1}p^{\alpha_{2}}_{2}...p^{\alpha_{q}}_{q}...p^{\alpha_{r}}_{r}...p^{\alpha_{s}}_{s}}{\left(p^{\alpha_{1}}_{1}p^{\alpha''_{2}}_{2}...p^{\alpha''_{q}}_{q}\right)\left(p^{\alpha'_{1}}_{1}p^{\alpha'_{2}}_{2}...p^{\alpha'_{r}}_{r}\right)}=\\
=\left(p^{\alpha_{1}-\alpha'_{1}-\alpha''_{1}}_{1}p^{\alpha_{2}-\alpha'_{2}-\alpha''_{2}}_{2}...p^{\alpha_{q}-\alpha'_{q}-\alpha''_{q}}_{q}...p^{\alpha_{r}-\alpha'_{r}}_{r}\right)...p^{\alpha_{s}}_{s}=	\\
=M\times O
\end{eqnarray*}
Where in $M$ there are all the primes in common and $O$ are the primes $p$ that divides $F_{n+2}$ but $p\not| F_{n+1}$ and $p\not| F_{n}$.  
So if $A'_{n+2}$ is the set of all the primes such $p \in A_{n+2}$ but $p \notin A_{n+1}$ and $gcd(F_{n+2},F_{n})=gcd(F_{n+2},F_{n+1})=gcd(F_{n+1},F_{n})=1$ then $\forall p \in A'_{n+2}$ $p | \binom{F_{n+2}}{F_{n}}$. $\qed$.




\bibliographystyle{spbasic} 

\end{document}